\documentclass[12pt,twoside]{amsart}
\usepackage{amssymb}

\nonstopmode

\numberwithin{equation}{section} 

\font\tengothic=eufm10 scaled\magstep 1 \font\sevengothic=eufm7
scaled\magstep 1
\newfam\gothicfam
      \textfont\gothicfam=\tengothic
      \scriptfont\gothicfam=\sevengothic

\newtheorem{theorem}{Theorem}[section]

\newtheorem{proposition}[theorem]{Proposition}
\newtheorem{corollary}[theorem]{Corollary}
\newtheorem{conjecture}[theorem]{Conjecture}

\theoremstyle{definition}
\newtheorem{remark}[theorem]{Remark}
\newtheorem{example}[theorem]{Example}

\newtheorem{question}[theorem]{Question}

\newcommand{\coker}{\operatorname{coker}}

\newcommand{\proj}[1]
{ \mathchoice
           { {\mathbb P}^{#1} }
           { {\mathbb P}^{#1} }
           { {\mathbb P}^{#1} }
           { {\mathbb P}^{#1} }
         }

\begin{document}
\title[Ideals of general forms and the ubiquity of the Weak Lefschetz
property]{Ideals of general forms and the ubiquity of the Weak Lefschetz
property}

\author[J.\ Migliore, R.M.\ Mir\'o-Roig]{J.\ Migliore$^*$, R.M.\
Mir\'o-Roig$^{**}$}
\address{Department of Mathematics,
        University of Notre Dame,
        Notre Dame, IN 46556,
        USA}
\email{Juan.C.Migliore.1@nd.edu}
\address{Facultat de Matem\`atiques,
Departament d'Algebra i Geometria, Gran Via de les Corts Catalanes
585, 08007 Barcelona, SPAIN } \email{miro@cerber.mat.ub.es}

\date{\today}
\thanks{$^*$ Partially supported by the University of Barcelona. \\
$^{**}$ Partially supported by BFM2001-3584.}



\begin{abstract}
Let $d_1,\dots,d_r$ be positive integers and let $I = (F_1,\dots,F_r)$ be an
ideal generated by forms of degrees $d_1,\dots,d_r$, respectively, in a
polynomial ring $R$ with $n$ variables.  With no further information virtually
nothing can be said about $I$, even if we add the assumption that $R/I$ is
Artinian.  Our first object of study is the case where the $F_i$ are chosen
generally, subject only to the degree condition.  When all the degrees are the
same we give a result that says, roughly, that they have as few first syzygies as
possible.  In the general case, the Hilbert function of $R/I$ has been
conjectured by Fr\"oberg.  In a previous work the authors showed that in many
situations the minimal free resolution of $R/I$ must have redundant terms which
are not forced by Koszul (first or higher) syzygies among the $F_i$ (and hence
could not be predicted from the Hilbert function), but the only examples came
when $r=n+1$.  Our second main set of results in this paper show that examples
can be obtained when $n+1 \leq r \leq 2n-2$.  Finally, we show that if
Fr\"oberg's conjecture on the Hilbert function is true then any such redundant
terms in the minimal free resolution must occur in the top two possible degrees
of the free module.
\par
Closely connected to the Fr\"oberg conjecture is the notion of Strong Lefschetz
property, and slightly less closely connected is the Weak Lefschetz property. 
We also study an intermediate notion, called the Maximal Rank property.  We
continue the description of the ubiquity of these properties, especially the
Weak Lefschetz property.  We show that any ideal of general forms in
$k[x_1,x_2,x_3,x_4]$ has the Weak Lefschetz property.  Then we show that for
certain choices of degrees, any complete intersection has the Weak Lefschetz
property and any almost complete intersection has the Weak Lefschetz property. 
Finally, we show that most of the time Artinian ``hypersurface sections'' of
zeroschemes have the Weak Lefschetz property.
\end{abstract}


\maketitle

\tableofcontents


 \section{Introduction} \label{intro}

Let $R = k[x_1,\dots,x_n]$ where $k$ is an infinite field.  Let $A
= R/I = \bigoplus_{i=0}^r A_i$ be a standard graded Artinian
$k$-algebra. The Weak Lefschetz property says that for a general
linear form $L$, the induced multiplication $(\times L) : A_i
\rightarrow A_{i+1}$ should have maximal rank, for each $i$.  The
Strong Lefschetz property says that for any power $d$, the
multiplication $(\times L^d) : A_i \rightarrow A_{i+d}$ has
maximal rank.  Notice that by semicontinuity the Strong Lefschetz property
 {\em implies} that for a  general form $F$ of arbitrary degree $d$, the induced
multiplication $(\times F): A_i \rightarrow A_{i+d}$ has maximal rank.  We will
call this latter property the {\em Maximal Rank property}. We do not know an
example where the Maximal Rank property is not equivalent to the Strong Lefschetz
property. 

It is well-known that the Weak Lefschetz property does not imply the Strong
Lefschetz property (see for instance \cite{HMNW}).  One would certainly 
expect, however, that  ``most'' Artinian $k$-algebras would have both
properties.  Many results in the last several years have contributed to
making this expectation more precise. On the other hand, many very natural
questions remain open.  We first recall several of these results and open
questions.

It was shown by R.\ Stanley \cite{stanley} and by J.\ Watanabe
\cite{watanabe1} that a monomial complete intersection always has
the Strong Lefschetz property.  By semicontinuity, it follows that
a {\em general} complete intersection has the property as well.  A
surprising step came in \cite{HMNW} where it was shown that for
$n=3$, {\em every} complete intersection has the {\em Weak
Lefschetz property}, extending a result of J.\ Watanabe
\cite{watanabe2}.  The problem remains open for $n \geq 4$, and
also for $n=3$ in the case of the Strong Lefschetz property.  It
was shown in \cite{HMNW} that any Artinian algebra in $k[x_1,x_2]$
has the Strong Lefschetz property.

Another interesting problem is to determine if the property holds
for Gorenstein Artinian $k$-algebras.  It was shown by J.\
Watanabe (\cite{watanabe1}, Example 3.9) that in any codimension,
``most'' Artinian Gorenstein rings possess the Strong  Lefschetz
property; more precisely, Watanabe showed that this holds for an open subset of
the projective space parameterizing the Artinian Gorenstein ideals
with fixed socle degree.  (Note that he does not show it for arbitrary Hilbert
function, and in fact the algebras that he produces are compressed, i.e.\ have
maximal Hilbert function.)  In
$k[x_1,x_2,x_3]$ it is not known if all Artinian Gorenstein ideals possess the
property, or if it at least holds for a general Artinian Gorenstein ideal with
fixed Hilbert function (cf.\ \cite{diesel}).  The same questions can
also be asked for Artinian Gorenstein ideals in
$k[x_1,\dots,x_n]$, possibly restricting to some subclass of such
ideals.  Watanabe \cite{watanabe3} proved a number of other strong
consequences of the Strong Lefschetz property for Gorenstein
rings.  On the other hand, it is known that not all Artinian
Gorenstein $k$-algebras have the Weak Lefschetz property, if $n
\geq 4$ (cf.\ for instance \cite{ikeda} Example 4.4).

Yet another interesting area concerns ideals of $r$ general forms,
with $r \geq n+1$ and the forms not necessarily of the same
degree.  Suppose that $I$ is such an ideal.  One could ask for the
Hilbert function of the algebra $R/I$ or even the minimal free
resolution of $R/I$.  It turns out that the Weak Lefschetz
property and especially the Maximal Rank property are
intimately connected to these questions.  One goal of this paper is to explore
these connections.

There are conjectures about the Hilbert function, due to R.\
Fr\"oberg, and about the minimal free resolution, due to A.\
Iarrobino, of these algebras.  The Fr\"oberg conjecture is
equivalent to the Maximal Rank property for such an algebra.  Anick
\cite{anick} settled the Hilbert function question for the case $R =
k[x_1,x_2,x_3]$, for any $r$, thus proving that {\em any} such $R/I$ has the
Maximal Rank property.  It is open in the case of more variables.

The Iarrobino conjecture said that the minimal free resolution of
$R/I$ should have no redundant terms apart from certain ones that
arise from Koszul syzygies.  This was disproved in a paper of the
present authors \cite{MMR2}, who analyzed the case $r = n+1$,
finding the explicit resolution in many cases and bounds in other
cases.  (At about the same time, counterexamples were found also by Pardue and
Richert \cite{pardue-richert}.)  In particular, in
\cite{MMR2} a connection was made to certain Gorenstein algebras (tying in with
the problem mentioned above) and a crucial step was the observation that these
algebras have the Strong Lefschetz property (although in this case it was enough
that they have the Weak Lefschetz property).

One of the few results on syzygies of ideals of general forms prior to
\cite{MMR2} was due to Hochster and Laksov \cite{hochster-laksov}.  It said that
an ideal of $r$ general forms of the same degree, $d$, spans a vector space of
maximum possible dimension in degree $d+1$.  In section 2 we extend this result
(Proposition \ref{gen hochster laksov}).  We are interested in the following
question: if $F_1,\dots,F_r$ are general forms of degree $d$ in
$k[x_1,\dots,x_n]$, what conditions force $R_t \cdot (F_1,\dots,F_r)$ to span a
vector space of maximal dimension in $(F_1,\dots,F_r)_{d+t}$?  We show that this
happens for any $t\le \min(d,t_0)$ where $t_0=\max\{l \ : \ 
\binom{d+l+2}{2}-(r-1)\binom{l+2}{2}\ge 0 \}$. In \cite{aubry}, M.\ Aubry
obtained a similar result, but relying on different assumptions.  Furthermore,
our proof is shorter and more elementary, and in certain ranges  of $n, r$
it improves his result. See Remark \ref{compare with aubry} for more details. 
Note that both results apply, in particular, only when $t<d$.

The third section contains our main results.  We are interested in the question
of trying to describe as well as possible the minimal free resolution of an
ideal of general forms.  In particular, when can they have redundant (``ghost'')
terms which are not related to Koszul syzygies (``non-Koszul ghost terms''), and
where can these ghost terms occur?  Since the only known non-Koszul ghost terms
occurred in the case of almost complete intersections ($n+1$ forms in $n$
variables), it was of interest to describe situations when more than
$n+1$ forms have non-Koszul ghost terms.  This is done in the first part of
Section 3, primarily with Theorem \ref{ghost terms in middle} and Corollary
\ref{ghost terms at end}.  In particular, we show that ghost terms can occur
(for the right choice of the degrees of the generators) when $n+1 \leq r \leq
2n-2$.

The next natural question is to narrow down where non-Koszul ghost terms can
possibly occur in an ideal of general forms.  All of the results mentioned so
far have the ghost terms occurring at the end of the resolution, or multiple
ghost terms occurring in a string of free modules starting at the end of the
resolution.  We show in Corollary \ref{not necessarily at the end} that this is
not always the case, although the examples created are somewhat special and
depend on generators of degree 2.

Of broad interest is the following question: in what degrees can the non-Koszul
ghost terms occur in a particular free module in the minimal free
resolution?  We give a conjecture that the syzygies can only be non-Koszul in
the top two degrees in each free module, and we prove this conjecture under the
hypothesis of Maximal Rank property (Conjecture \ref{conjecture on ghosts}, 
Proposition \ref{only ghost koszul in n var} and Proposition \ref{SLP implies
conj}).  Pardue and Richert \cite{pardue-richert} have recently obtained a
similar result, but the method of proof is entirely
different.  Chandler \cite{chandler} also has related work, in a very geometric
setting.

In particular, combining the results of \cite{MMR2}, the work mentioned above
and the result of Anick \cite{anick} that the Maximal Rank property holds
when $n=3$, we answer most of these questions for this case.

Section 4 contains results about the Weak Lefschetz property.  For example, one
could hope to prove that every Gorenstein ideal of height three possesses this
property by showing that the property is preserved under liaison. 
Unfortunately, we give a counterexample to this idea.  The main result of this
section is that every ideal of general forms in $k[x_1,x_2,x_3,x_4]$ has the
Weak Lefschetz property, again using Anick's result.

In Sections 5 and 6 we begin the task of seeing what ideals possess the Weak
Lefschetz property if we drop the assumption that they be ideals of general
forms.  A good first step was proved in \cite{HMNW}, where it was shown that
{\em every} complete intersection in $k[x_1,x_2,x_3]$ possesses this property. 
In Section 5 we show that under different assumptions (mostly on the degrees),
every complete intersection and every {\em almost complete intersection} has this
property.  The assumptions are quite restrictive, unfortunately.  Then in
Section 6 we show the Weak Lefschetz property for certain Artinian rings
obtained from zeroschemes.

The authors are grateful to Tony Iarrobino for pointing out the relevance of
Aubry's work to our Proposition \ref{gen hochster laksov}.


\section{A remark about first syzygies of general forms}

A problem that comes up surprisingly often in Algebra and Geometry
and which is closely related to the Strong Lefschetz property is
to determine the Hilbert series of the graded quotient $A=R/I$,
that is the series
\[Hilb_A(t)=\sum_{s=0}^\infty \dim_kA_st^s \]
\noindent where $A_s$ is the s-th graded piece of $A$. If $r\le n$
then $I$ is a complete intersection and the result is well known.
So, assume $r>n$, which in particular means that $A$ is Artinian.
In 1985, R.\ Fr\"oberg conjectured
\[
Hilb_A(t)= \left [ \frac{\prod _{i=1}^r(1-t^{d_{i}})}{(1-t)^n}
\right ]
\]
where $[\sum_{j=0}^{\infty}a_j t^{j}]=\sum_{j=0}^{\infty}b_j t^j$
with
\[
a_j  = \left \{
\begin{array}{ll}
  b_j & \mbox{ if } a_i\ge 0 \mbox{ for all } i\le j; \\
 0 & \mbox{ otherwise. }
\end{array}
\right.
\]

Several contributions to this apparently simple problem have been
made and there are at least three ways to attack this conjecture.
First, one could bound the number of variables. The conjecture was
proved to be true for $n=2$ in R.\ Fr\"oberg \cite{froberg} and
for $n=3$ in D.\ Anick \cite{anick}. Secondly, one could bound the
number of generators for the ideal $I$. The conjecture is easily
seen to be true for $r\le n$ and it was proved to be true for
$r=n+1$  by R.\ Stanley \cite{stanley}. It is also true if all the
generators have the same degree $d$ and $r\ge \frac{1}{n} {{d+n
\choose d+1}}$ (\cite{froberg} Example 4, p.\ 128). Thirdly, one
could prove that the conjecture is true for the first terms in the
Hilbert series. The first non-trivial statement comes for degree
$d+1$ with $d= \min \{ d_i \}$. In this degree the conjecture  is
equivalent to the following result of M.\ Hochster and D.\ Laksov:

\begin{proposition} Let $F_1,\dots , F_r$ be $r$ general forms of
degree $d$ in $R= k[x_1, \dots,x_n]$. Set $A=R/(F_1,\dots , F_r)$.
Then,
\[ dim_kA_{d+1}= \max \left \{0,{ n+d \choose d+1}-rn \right \} \]
\noindent  i.e., $\{ x_iF_j\}_{i=1,...,n;j=1,...,r }$ spans a
vector space of maximal dimension, namely,
\[
\min \left \{ rn,{n+d \choose d+1} \right  \}
\]
\end{proposition}

The goal of the next proposition is to extend the above result
about linear syzygies to higher degree syzygies.

\begin{proposition} \label{gen hochster laksov} Let $F_1,\dots , F_r$ be $r$
general forms of degree $d$ in $R= k[x_1, \dots,x_n]$. Set
$A=R/(F_1,\dots , F_r)$. Assume \[ {d+t_0+2 \choose
d+t_0}-(r-1){t_0+2 \choose t_0}\ge 0. \] \noindent Then,
\[ \dim_kA_{d+t}={n+d+t-1 \choose d+t}-r{t+n-1 \choose t}
\mbox{ for all } t\le t_0, \] i.e., $\{R_tF_j\}_{j=1,...,r}$ spans
a vector space of maximal dimension.
\end{proposition}

\begin{proof} We proceed by induction on $n$. For $n=2$ (resp.
$n=3$) the result is true and easily follows from Fr\"oberg
\cite{froberg} (resp. Anick \cite{anick}) and the fact that
$k[x_1,x_2]/(F_1,...,F_r)$ satisfies the Strong Lefschetz property (resp.
$k[x_1,x_2,x_3]/(F_1,...,F_r)$ satisfies the Maximal Rank property) .

Assume $n>3$. We want to construct $r$ forms of degree $d$, $
F_1,...,F_r$, such that
\[
\dim_kA_{d+t}={n+d+t-1 \choose d+t}-r{t+n-1 \choose t}
\]
being $A=k[x_1,\dots ,x_n]/(F_1,...,F_r)$. To this end, we first
consider $G_1,\dots , G_{r-1}$, a set of $r-1$ general forms of
degree $d$ in $k[x_1,\dots ,x_{n-1}] =: S$, and the ideal
$J=(G_1,\dots ,G_{r-1})\subset k[x_1,\dots ,x_{n-1}] = R$. By the
hypothesis of induction, for all $t\le t_0$, we have

\[
\dim_k (S/J)_{d+t} \ = \ {n+d+t-2 \choose d+t}-(r-1){t+n-2 \choose
t}.
\]

Now, we consider the ideal $I=(F_1, \dots ,F_r)\subset R $ where
\[
\begin{array}{rcl}
 F_i(x_1, \dots, x_n) & = & G_i(x_1, \dots ,x_{n-1})  \ \hbox{ for $i=1,\dots,
r-1,$ and } \\
 F_r(x_1, \dots, x_n) & = & x_n^d.
\end{array}
\]

We claim that for all $t\le t_0$, we have
\[
\begin{array}{rcl}
\displaystyle \dim_k (k[x_1,\dots ,x_{n}]/I)_{d+t} &  \le &
{n+d+t-1 \choose d+t} - \left [ (r-1){t+n-2 \choose t}+ {n+t-1
\choose n-1}+(r-1){t+n-2 \choose n-1} \right ] \\   \\ & = &
{n+d+t-1 \choose d+t}-r{t+n-1 \choose t}  .
 \end{array}
\]
To see this, consider the following subspaces of $R_{d+t}$:
\[
\begin{array}{rcl}
E_1 & = & S_t \cdot \langle F_1 ,\dots,F_{r-1} \rangle \\ E_2 & =
& R_t \cdot F_r  \\ E_3 & = & R_{t-1} \cdot \langle x_n F_1 ,
\dots, x_n F_{r-1} \rangle
\end{array}
\]
It is not hard to check that the three sets of canonical basis
elements are linearly independent, taken together.  This proves
the claim.  (Note that we have used here that  $t_0 < d$, so that
$E_3$ has no terms with a factor of $x_n^d$.  One can check that
this is true if $r \geq 5$, which holds here since we are assuming
$r>n>3$.)

Since we always have
\[
\dim_k (k[x_1,\dots ,x_{n}]/I)_{d+t} \ge {n+d+t-1 \choose
d+t}-r{t+n-1 \choose t}
\]
 we conclude that
\[
\dim_k (k[x_1,\dots ,x_{n}]/I)_{d+t} = {n+d+t-1 \choose
d+t}-r{t+n-1 \choose t}
\]
\end{proof}

\begin{remark} \label{compare with aubry}
M.\ Aubry has proved a result similar to Proposition \ref{gen
hochster laksov} (\cite{aubry}  Th\'eor\`eme 2.3) and it is worthwhile to make
some comments on the relation.  First, Aubry's result says that under certain
hypotheses the forms of degree $d$ span the maximum possible dimension in degree
$d+t$.  This could consist of the vector space of all forms of degree $d+t$, so
the obvious generators of $R_t \cdot I_d$ may span instead of being
linearly independent.  Our result gives different hypotheses to conclude only
that the forms are linearly independent.  

Fixing $t$, Aubry's result holds if
$d$ is larger than some function depending only on $n$, while ours depends only
on $r$.  The proof given above is shorter, and for some values of $n$ and $r$ it
improves on Aubry's result. 

For example, suppose that $n=10$ and we are interested in the span of $R_3\cdot
(F_1,\dots,F_r)_d$.  Aubry's result shows that the $r$ forms of degree
$d$ span the maximum dimension (independently of $r$) whenever 
\[
d \geq \frac{6(n-1)}{\sqrt[n-1]{(n-1)!}} - 3 + \frac{9}{\sqrt[n-2]{(n-2)!}} +
\frac{(n-1)^2}{\sqrt[n-1]{(n-1)!}} -n+5 \approx 27.
\]
If $d < 27$ his result does not apply (and indeed he remarks that it is not
the best bound possible).  If $n$ changes, the bound must be re-computed.

Our result above is independent of $n$, and says that the forms span the maximum
dimension whenever 
\[
\binom{d+5}{2} - (r-1)\binom{5}{2} \geq 0.
\]
We can thus choose any value of $d$ and the above inequality gives the values of
$r$ that allow us to reach our conclusion.  This range of $r$ works for any $n$.
\end{remark}


\section{Ghost terms in the minimal free resolution of an ideal of
general forms} \label{ghost term section}

Let $R= k[x_1, \dots,x_n]$ be a homogeneous polynomial ring over
an infinite field $k$ and let $I=(F_1, \dots , F_r)$ be an ideal
of $r$ generically chosen forms of degrees $d_i= \deg (F_i)$,
$i=1,\dots, r$.   We would like to comment on the minimal free
resolution and on the Hilbert function of such an ideal.

We begin with the minimal free resolution, which is the more
refined invariant of the two: knowing the minimal free resolution
gives the Hilbert function, but not conversely.  As with the
Hilbert function, for an ideal of $r$ general forms in $R =
k[x_1,\dots,x_n]$ there is an ``expected'' minimal free
resolution, conjectured  by Iarrobino \cite{iarrobino}.  This says
in effect that there should be no redundancies (``ghost terms'')
in the minimal free resolutions apart from those syzygies
(including higher syzygies) forced by Koszul relations among the
generators.  This was proven to be false in \cite{MMR2}, where it
was shown that in the case of $n=3$ and $r=4$ there can be
non-Koszul ghost terms.  Examples were given for larger values of
$n$, with $r = n+1$, but an examination of these examples shows
that the ghost terms appearing there arise (at least numerically)
from higher Koszul syzygies.

\begin{question} Is it the case that the only counterexample to Iarrobino's
conjecture comes when $n=3$ and $r=4$?
\end{question}

 We will show that this is not the case.  In fact, we can
find infinitely many counterexamples for any value of $n \geq 3$
and $n+1 \leq r \leq 2n-2$ (Corollary \ref{ghost terms at end}).
Another question is where the non-Koszul ghost terms can arise.
For example,

\begin{question}
Can there be non-Koszul ghost terms in the minimal free resolution
of ideals of general forms which do not arise between the last two
free modules in the resolution?
\end{question}

Our next result also answers this question in the affirmative.

\begin{theorem} \label{ghost terms in middle}
Let $R = k[x_1,\dots,x_n]$ and let $J = (F_1,\dots,F_n) \subset R$
be a complete intersection of general forms, with $\deg F_i = d_i$
for $1 \leq i \leq n$. Assume that $2 < d_1 \leq \dots \leq d_n$.
Let $d = d_1 + \dots + d_n$ and let $c = d - n - 1$.  Choose
general forms $F_{n+1}, \dots,F_{n+p}$ all of degree $c$, with $1
\leq p \leq n-2$.  Let $I = (F_1,\dots,F_n, F_{n+1}, \dots,
F_{n+p})$. Then for $j= p+1,\dots,n-1$,  $R/I$ has ghost terms
$R(-c-j)$ between the $j$-th and $(j+1)$-st free modules in the
resolution, which do not arise from any Koszul syzygies.
\end{theorem}

\begin{proof}

The last components of $R/J$ have dimension
\[
\dim \left ( R/J \right )_t = \left \{
\begin{array}{ll}
n & \hbox{ if $t = d-n-1$}; \\ 1 & \hbox{ if $t = d-n$}; \\ 0 &
\hbox{ if $t > d-n$}.
\end{array}
\right.
\]
Furthermore, $R/J$ has the Strong Lefschetz property thanks to
Stanley's result \cite{stanley}.  It follows that for the ideal
$I' = (F_1,\dots,F_n,F_{n+1})$, the last components of the ring
$R/I'$ have dimension
\[
\dim \left ( R/I' \right )_t = \left \{
\begin{array}{ll}
n-1 & \hbox{ if $t = d-n-1$}; \\ 0 & \hbox{ if $t \geq d-n$}.
\end{array}
\right.
\]
Then we get that
\[
\dim (R/I)_t = \left \{
\begin{array}{ll}
\dim (R/J)_t & \hbox{ if $t \leq d-n-2$}; \\ n-p & \hbox{ if $t =
d-n-1$}; \\ 0 & \hbox{ if $t \geq d-n$}.
\end{array}
\right.
\]

Let $G = [J:I]$ be the residual ideal.  By the Hilbert function
formula for linked Artinian rings (cf.\ \cite{DGO},
\cite{migliore}), the Hilbert function of $R/G$ is
\[
\dim (R/G)_t = \left \{
\begin{array}{ll}
1 & \hbox{ if $t = 0$}; \\ p & \hbox{ if $t = 1$}; \\ 0 & \hbox{
if $t > 1$}.
\end{array}
\right.
\]
It follows that the maximal socle degree of  $R/G$ is 1 and there
exist linear forms $L_1,...,L_{n-p}$ such that $(G)_1\cong
(L_1,...,L_{n-p})$. Hence,
\[
[Tor_{i}^R(G,k)]_{j} \cong [Tor_{i}^R((L_1,...,L_{n-p}),k)]_{j}
\]
for all $j\le 1+i-1=i$ and $R/G$ has a minimal free $R$-resolution
of the following type:
\begin{equation}\label{resol of G}
\begin{array}{ccccc}
0 \rightarrow R(-n-1)^{a_n} \rightarrow \dots \rightarrow
R(-n+p-2)^{a_{n-p+1}} \rightarrow
\begin{array}{ccccccc}
R(-n+p) \\ \oplus \\ R(-n+p-1)^{a_{n-p}}
\end{array}
\rightarrow \dots
\\ \\
\hbox{\hskip 2in} \dots \rightarrow
\begin{array}{ccccccccc}
R(-2)^{\binom{n-p}{2}} \\ \oplus \\ R(-3)^{a_2}
\end{array}
\rightarrow
\begin{array}{ccccccccc}
R(-1)^{n-p} \\ \oplus \\ R(-2)^{a_1}
\end{array}
\rightarrow R \rightarrow R/G \rightarrow 0.
\end{array}
\end{equation}
where $a_{i-1}$ is defined inductively by the equation
\[
\begin{array}{rcl}
0 & = & \dim (R/G)_i \\ & = & \displaystyle  {{n-1+i} \choose
{n-1}} +\sum_{j=1}^{i-1}(-1)^j \left [{{n-1+i-j} \choose {n-1}}
{{n-p} \choose j}+ a_j {{n-2+i-j} \choose {n-1}} \right ] \\ \\ &&
\displaystyle + \ (-1)^i {{n-p} \choose i}
\end{array}
\]
where we follow the convention that $\binom{a}{b} = 0$ if $a<b$,
so for example the last term is zero if $i > n-p$.

We get the diagram
\begin{equation} \label{setup for mapping cone}
\begin{array}{ccccccccccccccccccccccccc}
\dots \rightarrow & \displaystyle \bigoplus_{i<j} R(-d_i-d_j) &
\rightarrow & \displaystyle \bigoplus_{i=1}^n R(-d_i) &
\rightarrow & R & \rightarrow & R/J & \rightarrow  0 \\ &
\phantom{alpha_1} \downarrow \alpha_2 && \phantom{\alpha_1}
\downarrow \alpha_1 && \downarrow && \downarrow
\\
\dots \rightarrow & \displaystyle R(-2)^{\binom{n-p}{2}} \oplus
(\dots) & \rightarrow & \displaystyle R(-1)^{n-p} \oplus
R(-2)^{\binom{p+1}{2}} & \rightarrow & R & \rightarrow & R/G &
\rightarrow 0
\end{array}
\end{equation}
and the mapping cone construction gives a free resolution of
$R/I$. We get, after a small calculation,
\[
0 \rightarrow
\begin{array}{ccccccc}
R(-c-n)^\bullet \\ \oplus \\ R(-c-n+1)^\bullet
\end{array}
\rightarrow
\begin{array}{cccccccc}
R(-c-n+1)^\bullet \\ \oplus \\ R(-c-n+2)^\bullet
\end{array}
\rightarrow \dots \rightarrow
\begin{array}{cccccc}
R(-c-p-1)^\bullet \\ \oplus \\ R(-c-p)^\bullet
\end{array}
\rightarrow \dots \rightarrow R/I \rightarrow 0.
\]
The only chance for splitting off comes from redundancies induced
by the vertical maps in (\ref{setup for mapping cone}), and the
numerical assumption $d_1 > 2$ eliminates this possibility.  (See
for instance the proof of Corollary \ref{ghost terms at end}.)
\end{proof}

The ideal produced in Theorem \ref{ghost terms in middle} has a string of
ghost terms in the minimal free resolution. This string begins at the end and has
a length that depends on the number of generators.  We highlight the following
special case because it allows a simplification of the notation and it gives the
largest known (to us) number of generators of an ideal of general forms that has
a non-Koszul ghost term in the minimal free resolution.

\begin{corollary} \label{ghost terms at end}
Let $R = k[x_1,\dots,x_n]$ and let $J = (F_1,\dots,F_n) \subset R$
be a complete intersection of general forms, with $\deg F_i = d_i$
for $1 \leq i \leq n$. Assume that $2 < d_1 \leq \dots \leq d_n$.
Let $d = d_1 + \dots+ d_n$.  Choose general forms
$F_{n+1},\dots,F_{2n-2}$ all of degree $d-n-1$.  Let $I = (F_1,\dots,F_n,
F_{n+1}, \dots, F_{2n-2})$.  Then
$R/I$ has a ghost term of the form $R(-d+2)$ occurring in the last
and the penultimate free modules in the resolution, and this ghost
term does not arise from any Koszul syzygies.
\end{corollary}

\begin{proof}

We continue to use the notation of the last proof.  Since $p \leq
n-2$, $G$ has at least two generators of degree 1 and hence at
least one first syzygy term $R(-2)$.  Since  $p \geq 1$, $G$ has
at least one generator in degree 2.  (In fact, the number of
generators in degree 1 is $n-p$, the number in degree 2 is
$\binom{p+1}{2}$ and the number of first syzygies of degree 2 is
$\binom{n-p}{2}$.)

As above, the mapping cone gives a free resolution of $R/I$ that
ends
\[
\displaystyle 0 \rightarrow R(-d+1)^{n-p} \oplus
R(-d+2)^{\binom{p+1}{2}}  \rightarrow R(-d+2)^{\binom{n-p}{2}}
\oplus (\dots) \oplus \bigoplus_{i=1}^n R(d_i-d) \rightarrow \dots
\]
Because $d_1 > 2$, there is clearly no splitting off possible (no
component of the map $\alpha_1$ is an isomorphism) and we obtain
our ghost terms.  Similarly, because $d_1 > 2$, it is impossible
for $d-2$ to equal either the sum of $n-1$ of the $d_i$ or $n$ of
the $d_i$, so none of the ghost components arise from Koszul
syzygies.
\end{proof}

\begin{remark} \label{some splitting}
The assumption that $d_1 > 2$ in Corollary \ref{ghost terms at
end}  can be weakened substantially.  All we need is that the map
$\alpha_1$ does not pick out all the generators of $G$ of degree
2. For this to happen, it is enough that the number of generators
of $J$ of degree 2 be $< \binom{p+1}{2}$.  In particular, this is
guaranteed if $2n < p^2 +p$.  We are not sure to what extent
weakening this hypothesis affects ghost terms in the middle of the
resolution.
\end{remark}

\begin{example}
As remarked above, there are often still ghost terms when we allow
$d_1 = 2$, but fewer than expected because there is splitting off
in (\ref{setup for mapping cone}).  If the number of generators of
degree 2 is $\geq \binom{p+1}{2}$, it can happen that there are no
non-Koszul ghost terms.  For example, taking $n=4$ and choosing
general forms of degrees 2,2,2,4,5,5 we get a Betti diagram using
Macaulay \cite{macaulay} as follows:

\medskip

\begin{verbatim}
; total:      1     6    13    10     2 
; --------------------------------------
;     0:      1     -     -     -     - 
;     1:      -     3     -     -     - 
;     2:      -     -     3     -     - 
;     3:      -     1     -     1     - 
;     4:      -     2    10     8     - 
;     5:      -     -     -     1     2 
\end{verbatim}

\medskip

\noindent Here $p=2$ and we expect one ghost term at the end of
the resolution, but it is not there.  The term $R(-6)$ common to
the second and third modules in the resolution is Koszul, as is
the term $R(-4)$ common to the first and second modules.
\end{example}

The counterexample of D.\ Eisenbud and S.\ Popescu \cite{EP} to
Lorenzini's Minimal Resolution conjecture \cite{lorenzini} has a
ghost term which arises in the middle of the resolution, and
nowhere else.  One wonders if this can happen for ideals of
general forms.

\begin{question} \label{Q about ghost not at end}
Can an ideal of general forms have non-Koszul ghost terms which
occur in only one spot in the resolution, other than at the end?
\end{question}

In another direction, we have the following conjecture which gives a different
restriction on where the ghost terms can occur.  Some work in this direction has
been done by Pardue and Richert \cite{pardue-richert} and by Chandler
\cite{chandler}.

\begin{conjecture} \label{conjecture on ghosts}
Let $I \subset R = k[x_1,\dots,x_n]$ be generated by $r$ general
forms.  Let $c$ be the maximal socle degree of $R/I$ (i.e.\ the
last degree in which $R/I$ is non-zero).  Then the only possible
non-Koszul ghost terms correspond to copies of $R(-c-i)$ between
the $i$-th and the $(i+1)$-st free modules, for $i \geq 2$.  For $i=1$ there are
no  non-Koszul ghost terms.
\end{conjecture}

There are four situations in which we have some progress on these latter two
questions:

\begin{itemize}
\item When $n=3$ and $r=4$ a complete description of the possible minimal free
resolutions, and in particular of the possible ghost terms, was given in
\cite{MMR2} (cf.\ Remark \ref{case n=3, r=4} below).

\item  When $n=3$ and $r>4$ we can give a negative answer to Question \ref{Q
about ghost not at end} (cf.\ Remark \ref{ans question for n=3}) and prove
Conjecture
\ref{conjecture on ghosts} (cf.\ Remark \ref{prove conj for n=3}), using
\cite{anick}.

\item When $n= 4, 5$ or $6$  and some of the generators have degree 2 we can
modify the arguments above to give an affirmative answer to Question \ref{Q
about ghost not at end}.

\item When the Maximal Rank property holds we can prove Conjecture
\ref{conjecture on ghosts} (see Theorem~\ref{SLP implies conj}),
and in fact we prove something stronger.
\end{itemize}

\begin{remark} \label{case n=3, r=4}
For a precise description of all possible ghost terms when  $n=3$
and $r=4$ the reader can look at
 \cite{MMR2}, where there are many examples of general almost complete
intersection ideals in $k[x,y,z]$ generated by  homogeneous forms
of degree different from those described in Theorem \ref{ghost terms in middle}
and Corollary \ref{ghost terms at end},   with ghost terms in its minimal free
$R$-resolution (of course none of them violates Conjecture \ref{conjecture on
ghosts}). For instance, we consider an almost complete intersection ideal
$I\subset k[x,y,z]$ generated by  3 general  forms of degree 5 and one
general form of degree 7. The minimal resolution of $I$ is
\[0 \rightarrow R(-12)^4\oplus R(-11) \rightarrow R(-10)^7\oplus R(-11) 
\rightarrow
 R(-5)^3\oplus R(-7) \rightarrow I \rightarrow 0.\]
\end{remark}

Our next result shows that if the Maximal Rank property holds in $R$ then in
any case there can be no non-Koszul ghost terms at the beginning of the
resolution.  This is part of Conjecture \ref{conjecture on ghosts}.

\begin{proposition} \label{only ghost koszul in n var}
Let $R = k[x_1, \dots, x_n]$.  Let $I = (F_1,\dots,F_r)$ be an ideal
of $r \geq n$ general forms, and suppose that the minimal free
resolution of $I$ is
\[
0 \rightarrow {\mathbb F}_n \rightarrow \dots \rightarrow {\mathbb F}_2
\rightarrow {\mathbb F}_1 \rightarrow I \rightarrow 0.
\]
Assume that  the Maximal Rank property holds for all ideals of fewer than $
r$ general forms in $R$ (for example, this is known to hold for any $r$ if $n=3$
thanks to \cite{anick}). Then the only ghost terms that arise between ${\mathbb
F}_2$ and ${\mathbb F}_1$ come from Koszul relations.
\end{proposition}

\begin{proof}
Let us suppose that $\deg F_i = d_i$ for $1 \leq i \leq r$, and
that $d_1 \leq \dots \leq d_r$.   We proceed by induction on $r$.
If $r=n$ the result is obvious since $I$ is a complete
intersection.  Let $r = n+1$.  If there is a ghost term, it means
that $d_{n+1}$ is equal to the degree of a syzygy of $F_1,F_2,\dots,F_n$.
But these have only Koszul syzygies, so we are done. (The case $n=3$, $r
= 4$ in general was studied in \cite{MMR2}  section 4, where it
was shown that the only non-Koszul ghost terms that can occur are
at the end of the resolution, and a numerical analysis was done to
describe the shifts: it turns out to be $R(-c-2)$ where $c$ is the
maximal socle degree of $R/I$.  Note that this supports Conjecture
\ref{conjecture on ghosts}.)

Now assume that $r \geq n+2$.  Suppose that there is a non-Koszul
ghost term. This means that there is a syzygy
\[
A_1 F_1 + \dots + A_{r-2}F_{r-2} + A_{r-1}F_{r-1} = 0
\]
and $\deg F_r = \deg A_{r-1} + \deg F_{r-1}$. We distinguish two
cases:

\medskip

\noindent {\bf Case 1:} Suppose that $A_{r-1} \in
(F_1,\dots,F_{r-2})$.  Then writing $A_{r-1} = B_1 F_1 + \dots +
B_{r-2} F_{r-2}$, we have
\[
(A_1 + B_1 F_{r-2}) F_1 + \dots + (A_{r-2} + B_{r-2}F_{r-1})
F_{r-2} = 0.
\]
Hence the syzygy is actually a syzygy for $F_1,\dots,F_{r-2}$.
Hence this ghost term also appears in the minimal free resolution
of  the ideal $(F_1,\dots,F_{r-2},F_r)$, so by induction this
ghost term must arise from Koszul relations.

\medskip

\noindent {\bf Case 2:} Suppose that $A_{r-1} \notin
(F_1,\dots,F_{r-2})$. Then the image of $A_{r-1}$ in
$R/(F_1,\dots,F_{r-2})$ is a non-zero element annihilated by the
general form $F_{r-1}$.  That is,
\[
0 \neq \bar A_{r-1} \in \ker \left ( \times F_{r-1} : [
R/(F_1,\dots,F_{r-2}) ]_{\deg A_{r-1}} \rightarrow
[R/(F_1,\dots,F_{r-2})]_{\deg A_{r-1} + \deg F_{r-1}} \right ).
\]
But by our hypothesis, $R/(F_1,\dots,F_{r-2})$ has
the Maximal Rank property.  Therefore, it follows that the
above map is surjective. Consequently,
\[
[R/(F_1,\dots,F_{r-2},F_{r-1})]_{\deg A_{r-1}+\deg F_{r-1}} = 0.
\]
But  this degree is precisely the degree of $F_r$,
from which it follows that $F_r$ is not a minimal generator of $I
= (F_1,\dots,F_r)$.  This contradiction completes the proof.
\end{proof}

\begin{remark} Without the hypothesis ``general" Proposition \ref{only
ghost koszul in n var} turns out to be false. Indeed, if we consider
$I=(x^2,xy,xz,y^3,z^3)\subset k[x,y,z]$, the minimal resolution of
I is
\[0 \rightarrow R(-4)\oplus R(-7) \rightarrow R(-3)^3\oplus R(-4)^2 \oplus R(-6) 
\rightarrow
 R(-2)^3\oplus R(-3)^2 \rightarrow I \rightarrow 0.\]
\end{remark}

\begin{remark} \label{ans question for n=3}
As noted, Anick \cite{anick} has shown that an ideal $I$ of $r$ general forms in
$k[x,y,z]$ has the Maximal Rank property, so Proposition \ref{only ghost
koszul in n var} applies in this case.  This means that in the minimal free
resolution
\[
0 \rightarrow {\mathbb F}_3 \rightarrow {\mathbb F}_2 \rightarrow {\mathbb F}_1
\rightarrow I \rightarrow 0,
\]
the only possible non-Koszul ghost terms come between ${\mathbb F}_3$ and
${\mathbb F}_2$, answering Question~\ref{Q about ghost not at end}.
\end{remark}

Now we give a partial answer to Question \ref{Q about ghost not at end} for $n =
4,5$ or $6$, following our ideas in Theorem~\ref{ghost terms in middle},
Corollary
\ref{ghost terms at end} and Remark \ref{some splitting}.  Note that Remark
\ref{ans question for n=3} above precludes any hope of such a result when
$n=3$.

\begin{corollary} \label{not necessarily at the end}
Let $J \subset R = k[x_1,\dots,x_n]$, $n>3$, be an Artinian complete
intersection of general forms, with $\deg F_i = d_i$ for $1 \le i
\le n$.   Let $d = d_1 + \dots + d_n$.  Let $0 \neq \mu_2$ be the
number of generators of $J$ which have degree 2.  Choose general
forms $F_{n+1},\dots,F_{n+p}$, all of degree $d-n-1$ and let $I =
(F_1,\dots,F_n,F_{n+1}, \dots,F_{n+p})$.  If $n-p = 3$ and $\mu_2
\geq \binom{p+1}{2}$ then the only non-Koszul ghost terms are of
type $R(-c-n+2)$ between the $(n-2)$-nd and $(n-1)$-st free
modules in the resolution of $R/I$.
\end{corollary}

Note that the hypotheses of this corollary imply, in particular,
that $n \le 6$ since we have $\binom{p+1}{2} \leq \mu_2 \leq n = p+3$.

\begin{proof}
In Theorem \ref{ghost terms in middle}, the resolution (\ref{resol
of G}) shows that the ghost terms for $R/G$ can only come in the
first three modules, since $n-p=3$. When we link to $I$,  as noted
in Remark \ref{some splitting}, quadrics can split off.  The
hypothesis on $\mu_2$ guarantees that $J$ has more quadric
generators than $G$ does, so in (\ref{setup for mapping cone}) the
mapping cone removes all ghost terms at the end of the resolution
of $R/I$, leaving only one place where they remain, as claimed.
(Note that the vertical map $\alpha_2$ in (\ref{setup for mapping
cone}) does not split off any terms.)

One detail that should be checked is that all the (quadric)
generators of $G$ which numerically could be split off via
$\alpha_1$ in fact do get split off. This can be checked by
starting with a suitably general $\bar G$ with the desired Hilbert
function and choosing a complete intersection $\bar J \subset G$
beginning with $\mu_2$ general quadrics.  Then linking gives the
resulting $\bar I$ with the claimed splitting off, so the general
$I$ does as well, by semicontinuity.
\end{proof}

\begin{remark}
We still do not know if there can be ghost terms for $n=3$, $r > 4$
or for other values of $(n,r)$ than those described above.
It is conceivable that one could prove the conjecture for $r$
sufficiently large with respect to $n$, by proving the Strong Lefschetz (or just
Maximal Rank) property in this case. 
\end{remark}

Now we prove Conjecture \ref{conjecture on ghosts} (and in fact something
stronger) when the Maximal Rank property is known to hold.  Pardue and
Richert have a similar result, but  the method of proof is
completely different.

\begin{theorem} \label{SLP implies conj}
Let $I = (F_1,\dots,F_r) \subset R = k[x_1,\dots,x_n]$ be an ideal of generally
chosen forms of degrees $d_1,\dots,d_r$.  Assume that any ideal of 
$< r$ general forms in $R$ has the Maximal Rank property.  
Consider a minimal free resolution 
\[
0 \rightarrow {\mathbb F}_n \rightarrow {\mathbb F}_{n-1} \rightarrow \dots
\rightarrow {\mathbb F}_2 \rightarrow {\mathbb F}_1 \rightarrow I
 \rightarrow 0.
\]
Let $c$ be the maximal socle degree of $R/I$.  Then  the $i$-th
free module ${\mathbb F}_i$ has the form
\[
\begin{array}{rcl}
{\mathbb F}_1 & = & \displaystyle \bigoplus_{i=1}^r R(-d_i), \\

{\mathbb F}_i & = & R(-c-i)^\bullet \oplus R(-c-i+1)^\bullet \oplus {\mathbb
K}_i \ \hbox{ for $i \geq 2$},
\end{array}
\]
where 
 $R(-t)^\bullet$ refers to an unspecified (possibly zero)
number of copies of $R(-t)$ and ${\mathbb K}_i$ is the module of $i$-th
Koszul syzygies of degree $\leq c+i-2$.
In particular, Conjecture \ref{conjecture on ghosts} holds; that is, if $c$ is
the maximal socle degree of $R/I$ then the only possible non-Koszul ghost terms
correspond to copies of $R(-c-i)$ between the $i$-th and the $(i+1)$-st free
modules, for $i \geq 2$.  
\end{theorem}

\begin{proof}
Implicit in our hypotheses is the assumption that $F_1,\dots,F_r$ are all
minimal generators of $I$, so no $d_i$ is ``too large'' with respect to the
preceding degrees.  In particular, the form of ${\mathbb F}_1$ is clear.  It is
also clear from the socle degree that for $i \geq 2$ 
\[
{\mathbb F}_i  =  R(-c-i)^\bullet \oplus R(-c-i+1)^\bullet \oplus \dots
\]
Note that we know the value of $c$ because $(F_1,\dots,F_{r-1})$  has the
Maximal Rank property, by hypothesis, so we know the Hilbert function of
$R/I$.  We have seen in Proposition
\ref{only ghost koszul in n var} that there is no non-Koszul ghost term between
${\mathbb F}_2$ and
${\mathbb F}_1$.

We will proceed by induction on $r$.  When $r=n$, $I$ is a complete intersection
so all syzygies are Koszul, and the result is trivially true.  So now assume
that $r>n$.  Let $I' = (F_1,\dots,F_{r-1})$. Let $c'$ be the maximal socle degree
of
$R/I'$.  Note that $c' \geq c$.  Furthermore, $d_r \leq c'$ since otherwise
$F_r$ would not be a minimal generator of $I$.  Consider the map
\[
\times F_r : (R/I')_{t-d_r} \rightarrow (R/I')_t
\]
for $t \geq 0$.  The Maximal Rank property says that for the first values of
$t$ this map is injective, and then for the remaining values of $t$ it is
surjective.  In particular, the cokernel is zero whenever the map is
surjective.  But since the cokernel is precisely $(R/I)_t$, we have that the map
is injective for $t \leq c$ and surjective for $t \geq c+1$.  

Now consider a syzygy of the generators of $I$, which we will write as follows:
\[
A_r F_r = A_1 F_1 + \dots + A_{r-1} F_{r-1}.
\]
If $\deg A_r + \deg F_r \leq c$ then injectivity forces $A_r \in I'$.  Then an
argument similar to that given in Proposition \ref{only ghost koszul in n var},
Case 1, shows that in fact the above syzygy can be written as 
\[
0 = (A_1 -B_1F_r)F_1 + \dots + (A_{r-1}-B_{r-1}F_r)F_{r-1}.
\]
Since $c \leq c'$, the inductive hypothesis shows that this is a Koszul syzygy. 
It follows that the only non-Koszul syzygies in fact correspond to copies of
$R(-c-1)$ and $R(-c-2)$, i.e.\ we have
\[
{\mathbb F}_2 = R(-c-2)^\bullet \oplus R(-c-1)^\bullet \oplus {\mathbb K}_2
\]
where ${\mathbb K}_2$ are only Koszul syzygies.  

Now consider ${\mathbb F}_3$ and suppose it has a component $R(-t)$ where $t
\leq c+1$.  Let $M_1$ be the module of first syzygies, so
\[
\dots \rightarrow 
\begin{array}{ccccccccc}
R(-c-3)^\bullet \\
\oplus \\
R(-c-2)^\bullet \\
\oplus \\
R(-t)^\bullet \\
\oplus \\
\vdots
\end{array}
\rightarrow
\begin{array}{ccccc}
R(-c-2)^\bullet \\
\oplus \\
R(-c-1)^\bullet \\
\oplus \\
{\mathbb K}_2
\end{array}
\rightarrow M_1 \rightarrow 0.
\]
Then any copy of $R(-t)$ is a syzygy of generators of $M_1$ corresponding to
summands of ${\mathbb K}_2$, i.e.\ is a Koszul second syzygy.  A similar
argument for the remaining free modules ${\mathbb F}_i$ completes the proof.
\end{proof}

\begin{remark} \label{prove conj for n=3}
Since Anick \cite{anick} has shown that the Maximal Rank property holds when
$n=3$, we have proven Conjecture \ref{conjecture on ghosts} for this case.
\end{remark}


\section{Some observations on the Weak Lefschetz property}

In this section we collect some general remarks.  First, in the
introduction it was asked whether all Gorenstein $k$-algebras in
$k[x_1,x_2,x_3]$ have the Weak (or Strong) Lefschetz property, as
was recently shown for complete intersections \cite{HMNW}.  A
natural way that one might hope to prove this result is by
liaison.  If one could show that the Weak Lefschetz property is
preserved under liaison, then the result of \cite{HMNW} and the
desired result for Gorenstein $k$-algebras would follow trivially
(since it was shown by Watanabe \cite{watanabe4} that a Gorenstein
ideal is in the liaison class of a complete intersection).

Unfortunately, it is not true that the Weak Lefschetz property is
preserved under liaison, as the following example shows.

\begin{example} \label{WLP not preserved under liaison}
Let $R = k[x_1,\dots,x_n]$ and let $I_1 = (x_1^2, x_1x_2,x_1x_3,
\dots, x_1x_n, x_2^3, x_3^3, \dots, x_n^3)$. Note that $R/I_1$
does not have the Weak Lefschetz property since $x_1 \in R/I_1$ is
annihilated by all linear forms. On the other hand, we claim that
$I_1$ is linked via the complete intersection $J_1 = (x_1^2,x_2^3,
\dots, x_n^3)$ to the ideal $I_2 = (x_1, x_2^3, x_3^3, \dots,
x_n^3, x_2^2x_3^2\cdots x_n^2)$, which in turn is linked via the
complete intersection $J_2 = (x_1, x_2^3, x_3^3, \dots, x_n^3)$ to
the ideal $I_3 = (x_1,x_2,x_3, \dots, x_n)$.

For the first link, the inclusion $I_2 \subset [J_1 :I_1]$ is
clear.  For the reverse inclusion, note first that $[J_1 :I_1]$ is
a monomial ideal since both $J_1$ and $I_1$ are monomial ideals.
Let $f = x_1^{a_1}x_2^{a_2}\cdots x_n^{a_n} \in [J_1:I_1]$.  We
want to show that $f \in I_2$.  Without loss of generality we may
assume that $a_1 =0$, $a_2 \leq 2,\dots,a_n \leq 2$ since
otherwise it is clear that $f \in I_2$.  But we have that $f \cdot
x_1x_i \in J_1$ for all $1 \leq i \leq n$.  From this, and our
assumption, it follows easily that $a_2 = \dots =a_n = 2$, so $f
\in I_2$.  The second link is left to the reader.
\end{example}

This example also serves to suggest the following.  Note that in
\cite{HMNW} it was shown that {\em every} Artinian ideal in
$k[x_1,x_2]$ has the Weak Lefschetz property.

\begin{question}
For any integer $n \geq 3$, find the  maximum number $A(n)$ (if it
exists) such that {\em every} Artinian ideal $I \subset
k[x_1,\dots,x_n]$ with $\mu (I) \leq A(n)$ has the Weak Lefschetz
property (where $\mu(I)$ is the minimum number of generators of
$I$).  Included in this question is whether every complete
intersection in $k[x_1,\dots,x_n]$ has the Weak Lefschetz
property, which would say that $A(n)$ exists and is $\geq n$.
\end{question}

\noindent Note that in \cite{HMNW} it was shown that every
complete intersection in $k[x_1,x_2,x_3]$ has the Weak Lefschetz
property, so $A(3) \geq 3$ (and in particular $A(3)$ exists).
Example \ref{WLP not preserved under liaison} shows that $A(n)
\leq 2n-2$ for any $n \geq 3$, if it exists.  We wonder if it is
true that $A(n) = 2n-2$.  In any case, the two most interesting
cases for now are to determine if every complete intersection in
$k[x_1,\dots, x_n]$ ($n \geq 4$) has the Weak Lefschetz property,
and if every almost complete intersection in $k[x_1,x_2,x_3]$ has
the Weak Lefschetz property. (This would say that $A(3) = 4$.  We
believe both of these to be true.  The results in this section and
(especially) the next are intended to contribute to the solution
of these questions.  For instance, Proposition \ref{WLP if last
gen is big} proves that every complete intersection in
$k[x_1,\dots,x_n]$ has the Weak Lefschetz property if the last
generator is of sufficiently large degree.

\vskip 2mm We would also like to remark on an easy consequence of
the previously-mentioned theorem of Anick \cite{anick}, that in the
ring $S = k[x_1,x_2,x_3]$, if $I$ is any ideal of general forms in $S$, then
$S/I$ has the Maximal Rank property.  Note that although we
state this result only for $k[x_1,x_2,x_3,x_4]$, the proof also holds for any
number of variables if the Maximal Rank property holds in a ring of one
fewer variables, a hypothesis similar to that used for instance in Proposition
\ref{only ghost koszul in n var} and Theorem \ref{SLP implies conj}. 

\begin{proposition} \label{anick implies WLP}
Any ideal of general forms in the ring $R=k[x_1,x_2,x_3,x_4]$ has
the Weak Lefschetz property.
\end{proposition}

\begin{proof}
Let $L \in R_1$ be a general linear form and let $S = R/(L) \cong
k[x_1,x_2,x_3]$.  Let $I = (F_1,\dots,F_{r-1}, F_r) \subset R$ be
an ideal of general forms, and write $I = I' + (F_r)$, where $I' =
(F_1,\dots,F_{r-1})$.  Suppose that $\deg F_i = d_i$ for $1 \leq i
\leq r$.  For $F \in R$ we denote by $\bar F$ the
restriction to $S$, and similarly for an ideal $J \subset R$.

The proof will be by induction on $r$.  The result is well known if $r = 4$, so we can assume
that $r \geq  5$.    Consider the diagram
\[
\begin{array}{cccccccccccccccccccc}
(R/I')_t & \stackrel{\times F_r}{\longrightarrow} & (R/I')_{t+d_r}
& \rightarrow & (R/I)_{t+d_r} & \rightarrow & 0 \\ \phantom{\times
L} \downarrow \times L && \phantom{\times L} \downarrow \times L
&& \phantom{\times L} \downarrow \times L  \\ (R/I')_{t+1} &
\stackrel{\times F_r}{\longrightarrow} & (R/I')_{t+d_r+1} &
\rightarrow & (R/I)_{t+d_r+1} & \rightarrow & 0 \\ \downarrow &&
\downarrow \\ (S/\bar I')_{t+1} & \stackrel{\times \bar
F_r}{\longrightarrow} & (S/\bar I')_{t+d_r+1} \\ \downarrow &&
\downarrow \\ 0 && 0
\end{array}
\]
Note that we do not assume that $\times F_r$ has maximal rank.
By induction, we may assume that $R/I'$ has the Weak Lefschetz property, i.e.\
that the first two vertical maps
$\times L$ have maximal rank (injective or surjective depending on
$t$).  If $t$ is such that the second vertical map $\times L$ is
surjective then $(S/\bar I')_{t+d_r+1} = 0$ and it is not hard to
see that the last vertical map $\times L$ is also surjective.  So
suppose that the second vertical map $\times L$ is injective.  By the Weak
Lefschetz property, it follows that the first vertical map $\times L$ is also
injective.  Then it
is tedious but not hard to show that if the last horizontal map $\times
\bar F_r$ is injective then the last vertical map $\times L$ is injective,
and if the last horizontal map $\times \bar F_r$ is surjective
then so is the last vertical map $\times L$.  Since by Anick's result
the last horizontal map $\times \bar F_r$ is always either
injective or surjective, we are done.
\end{proof}

\begin{remark} \label{force WLP}
Let $R = k[x_1,\dots,x_n]$ and let $R/I$ be an Artinian Gorenstein
ring. Without loss of generality assume that no generator of $I$
has degree 1 and that the Hilbert function of $R/I$ is
\[
1 \ \ n \ \ h_2 \ \ h_3 \ \  \dots \ \ h_3 \ \ h_2 \ \ n \ \ 1.
\]
Assume that the socle degree is $s$.  Then there are some
situations which guarantee that $R/I$ has the Weak Lefschetz
property:
\begin{itemize}

\item[(a)] $s$ is even and $h_i = \binom{n-1+i}{i}$ for $0 \leq i \leq
\frac{s}{2}$.  (This means that $R/I$ agrees with $R$ through the
first half of the Hilbert function; that is, $R/I$ is {\em compressed} with even
socle degree.)  Note that this does not hold if
$s$ is odd. Indeed, H.\ Ikeda  \cite{ikeda} has found an example of a
Gorenstein Artinian ring with Hilbert function $1 \ 4 \ 10 \ 10 \
4 \ 1$ which does not have the Weak Lefschetz property.  Of course
it fails ``in the middle."

\item[(b)]  $n=3$ and the Hilbert functions contains a sequence $t,t,t$ (at
least three) in the middle.  This is an easy consequence of
\cite{IK} Theorem 5.77 (a).  Note that for this result we do not
need the ``growth like $R$" assumed in part (a) above.

\item[(c)] $n=3$ and the skew-symmetric Buchsbaum-Eisenbud matrix \cite{BE} has
only linear entries.  It is possible to make a direct argument, but in
fact this is {\em equivalent} to the statement in part (a).  This
can be seen using Diesel \cite{diesel} or more simply using
Corollary 8.14 of \cite{MN3} (using the case $s=2t$).  This latter
result shows that even for larger $n$, having $s$ even and $h_i$
maximal guarantees a resolution of the form
\[
0 \rightarrow R(-s-n) \rightarrow R(-t-n+1)^{\beta_{n-1}}
\rightarrow \dots \rightarrow R(-t-1)^{\beta_1} \rightarrow R
\rightarrow R/I \rightarrow 0
\]
which is linear except at the beginning and end.
\end{itemize}

We would like to ask whether the condition in (a) also forces
$R/I$ to have the Strong Lefschetz property, or at least the Maximal Rank
property.
\end{remark}

\begin{remark}
Every height $n$ Artinian ideal $I \subset k[x_1,\dots,x_n]$ with
a linear resolution has the Strong Lefschetz property.  Indeed,
suppose the resolution has the form
\[
0 \rightarrow R(-p-n+1)^{a_n} \rightarrow \cdots \rightarrow
R(-p)^{a_1} \rightarrow R \rightarrow R/I \rightarrow 0.
\]
Then the socle degree of $R/I$ is $p-1$, so the Hilbert function
is
\[
h_{R/I}(t) = \left \{
\begin{array}{ll}
\binom{n-1+t}{n-1} & \hbox{if $t < p$} \\ 0 & \hbox{if $t \geq p$}
\end{array}
\right.
\]
That is, $R/I$ agrees with $R$ until degree $p-1$ and then is
zero, so the Strong Lefschetz property is clear.

One of the most basic open problems at this stage is whether every height three
Gorenstein ideal in $k[x_1,x_2,x_3]$ has the Weak Lefschetz property (or better,
the Strong Lefschetz property).  It is known that every height three
complete intersection has the Weak Lefschetz property (see Theorem \ref{HMNW
result}).  In Remark \ref{force WLP} we saw that it also holds for a height
three Gorenstein ideal with only linear entries in the Buchsbaum-Eisenbud
matrix.  We propose as the next step to prove that a height three Gorenstein
ideal with only quadratic entries in the Buchsbaum-Eisenbud matrix has the Weak
Lefschetz property.  A first example could be a Gorenstein ideal with minimal
free resolution
\[
0 \rightarrow R(-10) \rightarrow R(-6)^5 \rightarrow R(-4)^5 \rightarrow R
\rightarrow R/I \rightarrow 0.
\]
\end{remark}


\section{Almost complete intersections and the Weak Lefschetz property}

In \cite{MMR2} the authors considered ideals in $k[x_1,\dots,x_n]$
that have $n+1$ generators, chosen generically.  The goal was to
describe the minimal free resolution of such an ideal, and along
the way to describe ghost terms that arise (as we have generalized
in Section \ref{ghost term section} above).  

Note that such an ideal is an almost complete intersection.  In
this section we would like to explore to what extent an Artinian
almost complete intersection, whose generators are not necessarily
chosen generically, must have the Weak Lefschetz property.  To
begin, however, we consider complete intersections.

Let $R' = k[x_1 ,x_2 ,x_3]$ and consider a complete intersection
$I' = (F_1,F_2,F_3)$ whose generators have degrees $d_1 \leq d_2
\leq d_3$.  In \cite{watanabe2} Corollary 2, Watanabe has shown
that if $d_3 \geq d_1 + d_2 -3$ then $R'/I'$ has the Weak
Lefschetz property.  This was generalized as follows:

\begin{theorem}[\cite{HMNW}] \label{HMNW result}
Every complete intersection in $R'$ has the Weak Lefschetz
property.
\end{theorem}

It is an open problem to show that every complete intersection $I
= (F_1,\dots,F_n)  \subset R = k[x_1,\dots,x_n]$ has the Weak
Lefschetz property. However, we would like to remark that at least
Watanabe's result extends to $R$ (in a slightly weaker form).

\begin{proposition} \label{WLP if last gen is big}
Let $I = (F_1 ,\dots,F_{n-1},F_n) \subset R = k[x_1,\dots,x_n]$
be a complete intersection.   Suppose that $d_i = \deg F_i$ for $1
\leq i \leq n$, and $2 \leq d_1 \leq \dots \leq d_n$.  Assume that
one of the following holds:
\begin{itemize}
\item[a.] $d_1 + \dots +d_{n-1} +d_n -n$ is even and
\[
d_n > d_1 + \dots + d_{n-1} - n;
\]

\item[b.] $d_1 + \dots + d_{n-1} + d_n -n$ is odd and
\[
d_n > d_1 + \dots + d_{n-1} - n +1.
\]
\end{itemize}
Then  $R/I$ has the Weak Lefschetz property.
\end{proposition}

\begin{proof}
Throughout this proof we set $J = (F_1,\dots,F_{n-1})$.  We denote
by $\bar R$ the ring $R/(L)$ for a general linear form $L$, and by
$\bar F$ (resp.\ $\bar J$) the restriction to $\bar R$ of a
homogeneous polynomial $F$ (resp.\ the homogeneous ideal $J$).
Note that $J$ is the ideal of a zeroscheme $Z$ in $\proj{n-1}$ (a
complete intersection), and that  $h_{R/J}(t) = \deg Z$ for $t
\geq d_1 + \dots + d_{n-1} -(n-1)$, since this is the socle degree
of $\bar R/ \bar J$.  In any case, $(R/J)_{t-1} \rightarrow
(R/J)_t$ is injective for all $t$, since $R/J$ is the coordinate
ring of a zeroscheme.  Note also that the socle degree of $R/I$ is
$d_1+\dots+d_{n-1} + d_n -n$.

Let us first assume that $d_1 + \dots + d_{n-1} + d_n -n$ is even
and $d_n > d_1 + \dots + d_{n-1} - n$. Then the Hilbert function
$h_{R/I}$ is symmetric, and the midpoint is in degree $
\frac{d_1+\dots+d_{n-1}+d_n-n}{2}$.  Note that the hypothesis $d_n
> d_1+\dots+d_{n-1} -n$ is equivalent to $
\frac{d_1+\dots+d_{n-1}+ d_n -n}{2} < d_n$.

Now, because $(R/I)_t = (R/J)_t$ for $t < d_n$, the multiplication
$(R/I)_{t-1} \rightarrow (R/I)_t$ induced by a general linear form
is injective for $t < d_n$.  Since the midpoint occurs in degree $
\frac{d_1+\dots+d_{n-1}+d_n-n}{2} < d_n$, we have injectivity in
``the first half.''  By duality, this is enough to prove the Weak
Lefschetz property for $R/I$.

If $d_1 + \dots + d_{n-1} + d_n -n$ is odd and $d_n > d_1 + \dots
+ d_{n-1} - n +1$, there is a small additional problem to
overcome.  In this case, the midpoint of $h_{R/I}$ is not an
integer, and we have to prove injectivity until just past the
midpoint.  That is, we have to show injectivity of $(R/I)_{t-1}
\rightarrow (R/I)_t$ for $t \leq
\frac{d_1+\dots+d_{n-1}+d_n-n+1}{2}$. If $d_n > d_1 + \dots
+d_{n-1}-n+1$, the same argument as the even case gives the
result.
\end{proof}

We now turn to almost complete intersections.  We first consider
the ring $R = k[x_1,x_2,x_3]$.  For any real number $x$, we set
$\lceil x \rceil = \min \{ n \in {\mathbb Z} \ | \ n \geq x \}$.

Let $I = (F_1,F_2,F_3,F_4)$ be an Artinian almost complete
intersection in $R$.  Note that while there is no loss of
generality in assuming that three of the four generators form a
regular sequence, say $(F_1,F_2,F_3)$, it  {\em does} become a
restriction if we further impose the condition $d_1 \leq d_2 \leq
d_3 \leq d_4$, since that forces us to assume that the three
generators of least degree form a regular sequence.  So {\em we do
{\bf not} make the restriction $d_1 \leq d_2 \leq d_3 \leq d_4$}.

\begin{proposition} \label{lambda minus 1}
Let $I = (F_1,F_2,F_3,F_4) \subset R$ be a height three almost
complete intersection Artinian ideal with generators of degrees
$d_1,d_2,d_3,d_4$. Assume that  $d_4 \geq \lceil
\frac{d_1+d_2+d_3}{2} \rceil -1$ and that $(F_1,F_2,F_3)$ form a
regular sequence.  Then $R/I$ has the Weak Lefschetz property.
\end{proposition}

\begin{proof}
Let $J = (F_1,F_2,F_3) \subset R$.  Note that the socle degree of
$R/J$ is $d_1+d_2+d_3-3$.   If $d_4 \geq d_1+d_2+d_3-2$ then $I=J$
so the result follows from Theorem \ref{HMNW result}.  So we
assume that
\[
\lambda -1 := \left \lceil \frac{d_1+d_2+d_3}{2} \right \rceil -1
\leq d_4 < d_1+d_2+d_3-2.
\]
 The hypothesis on
$d_4$, together with Theorem \ref{HMNW result}, show that for a
general linear form $L$ and for any $t \leq \lambda-2$, we have an
injection
\[
(R/I)_{t-1} = (R/J)_{t-1} \stackrel{\times L}{\longrightarrow}
(R/J)_t = (R/I)_t.
\]
On the other hand,  $\times L : (R/J)_{t-1} \rightarrow (R/J)_{t}$
is surjective  for all $t \geq \lambda-1$ since $R/J$ has the Weak
Lefschetz property, by Theorem \ref{HMNW result}.  But we also
have a surjection $(R/J)_t \rightarrow (R/I)_t$ for all $t$. Then
for $t \geq \lambda-1$ we have the commutative diagram
\[
\begin{array}{ccccccccccccc}
(R/J)_{t-1} & \rightarrow & (R/J)_t & \rightarrow & 0 \\
\downarrow && \downarrow \\ (R/I)_{t-1} & \rightarrow & (R/I)_t &
\\ \downarrow && \downarrow \\ 0 && 0
\end{array}
\]
from which the surjectivity of $\times L : (R/I)_{t-1} \rightarrow
(R/I)_t$ follows immediately.
\end{proof}

Proposition \ref{lambda minus 1} can be generalized to $n$ variables and we have

\begin{proposition}Let $I = (F_1,F_2,...,F_{n+1}) \subset k[x_1,...,x_n]$
 be a height $n$ almost
complete intersection Artinian ideal with generators of degrees
$d_1,...,d_n,d_{n+1}$. Assume that  $d_{n+1} \geq \lceil
\frac{d_1+...+d_n}{2} \rceil -1$,  that $(F_1,...,F_n)$ form a
regular sequence and that $d_n > d_1 + \dots + d_{n-1} - n$ (resp.
$ d_n > d_1 + \dots + d_{n-1} - n +1$) if $d_1 + \dots +d_{n-1}
+d_n -n$ is even (resp. $d_1+\dots+d_{n-1}+d_n -n$ is odd).  Then
$k[x_1,...,x_n]/I$ has the Weak Lefschetz property.
\end{proposition}

\begin{proof} It is analogous to the proof of Proposition
\ref{lambda minus 1} using Proposition \ref{WLP if last gen is big} instead of
Theorem \ref{HMNW result}.
\end{proof}

\begin{proposition} \label{lambda minus 2}
Let $I = (F_1,F_2,F_3,F_4) \subset k[x_1,x_2,x_3] = R$ be an Artinian almost
complete intersection and assume that $J = (F_1,F_2,F_3)$ forms a regular
sequence.   Suppose that $2 \leq d_i$ for each $i = 1, 2, 3, 4$ and that 
\[
d_4 = \left \lceil \frac{d_1+d_2+d_3}{2} \right \rceil -2 := \lambda-2.
\]
Finally, suppose that there exists a linear form $L$ such that 
\[
F_4 \notin \ker
\left [ (R/J)_{\lambda-2} \stackrel{\rho}{\longrightarrow} (\bar R/\bar
J)_{\lambda-2} \right ],
\]
 where $\rho$ is the restriction modulo $L$.  Then
$R/I$ has the Weak Lefschetz property.
\end{proposition}

\begin{proof}
We remark that we believe that the last hypothesis always holds, and hence is
superfluous, but we are not able to prove this.

Note that the Hilbert function of $R/J$ is symmetric.  If $d_1+d_2+d_3$ is odd
then there is a well-defined middle term in degree $\lambda-2$, and by Theorem
\ref{HMNW result} for a general linear form $L$ we have an injection
$\times L : (R/J)_{i-1} \rightarrow (R/J)_i$ for all $i \leq \lambda-2$ and
surjection for all $i \geq \lambda-1$.  If 
$d_1+d_2+d_3$ is even then there are (at least) two equal terms in the middle,
in degrees $\lambda-2$ and $\lambda-1$, and again we have an injection
$\times L : (R/J)_{i-1} \rightarrow (R/J)_i$ for all $i \leq \lambda-2$, and
now also an isomorphism $\times L : (R/J)_{\lambda-2} \rightarrow
(R/J)_{\lambda-1}$.  Note that $d_4 = \lambda-2$.

Arguing as in Proposition~\ref{lambda minus 1}, to complete the proof it is
enough to see that for a general linear form $L$, the induced map
\begin{equation}\label{needed maxrk}
(R/I)_{d_4-1} \stackrel{\times L}{\longrightarrow} (R/I)_{d_4}, \ \ \
\hbox{i.e.} \
\ \ (R/I)_{\lambda-3} \stackrel{\times L}{\longrightarrow} (R/I)_{\lambda-2},
\end{equation}
has maximal rank.  

\medskip

\noindent {\bf Case 1.} Suppose that $d_1+d_2+d_3$ is even and $d_3 > \lambda$,
or equivalently that $d_3 > d_1+d_2$.  Then
\[
\begin{array}{rcl}
d_1+d_2 & = & d_1+d_2+d_3 - d_3 \\
& = & 2 \lambda - d_3 \\
& < & \lambda.
\end{array}
\]
  Note that we have a Hilbert function 
\[
h_{R/(F_1,F_2)}(t) = d_1 d_2 \ \ \hbox{for} \ \  t \geq d_1+d_2-2
\]
and that $(F_1,F_2)$ is the saturated ideal of a zeroscheme in $\proj{2}$. 
Therefore the Hilbert function of $R/J$ has terms

\begin{center}
\begin{tabular}{c|cccccccccccccccccccccc}
$t$ & $\dots$ & $\lambda-3$ & $\lambda -2$ & $\lambda-1$ & $\lambda$ \\
\hline
$h_{R/J}(t)$ & $\dots$ & $d_1d_2$ & $d_1d_2$ & $d_1d_2$ & $d_1d_2$ & $\dots$
\end{tabular}
\end{center}
Hence by Theorem \ref{HMNW result} we have (in particular) a
surjection $\times L : (R/J)_{\lambda-3} \rightarrow (R/J)_{\lambda-2}$, so the
same proof as in Proposition \ref{lambda minus 1} gives the surjection
(\ref{needed maxrk}).

\medskip

\noindent {\bf Case 2.}  
Suppose that $d_1+d_2+d_3$ is odd and $d_3 > \lambda -1$,
or equivalently that $d_3 > d_1+d_2-1$.  Then as in Case 1, we quickly check
that $d_1+d_2 < \lambda$.  Now the Hilbert function calculation of Case 1 is
the same in degrees $\lambda-3, \lambda-2$ and $\lambda-1$ (but could change in
degree $\lambda$ if $d_3 = \lambda$).  But then the proof is identical to that
of Case 1. 

\bigskip

When we begin with $(F_1,F_2, F_3)$ and add the generator $F_4$ in degree
$\lambda-2$,  the Hilbert function of $R/J$ is unchanged in degrees $\leq
\lambda-3$ and drops by 1 in degree $\lambda-2$.  Cases 1 and 2 cover the only
situations where $h_{R/J}(\lambda-3) \geq h_{R/J}(\lambda-2)$ (in fact it is
$=$).  In all other cases $h_{R/J}(\lambda-3) < h_{R/J}(\lambda-2)$, so it is
enough to show for (\ref{needed maxrk}) that we have an injection $\times L :
(R/I)_{\lambda-3} \rightarrow (R/I)_{\lambda-2}$.  Note that we have the
corresponding injection for $R/J$ by Theorem \ref{HMNW result}.

To this end, we consider the commutative diagram
\begin{equation}\label{bigdiagram}
\begin{array}{ccccccccccccccc}
&& 0 && 0  \\
&& \downarrow && \downarrow \\
0 & \rightarrow & E & \rightarrow & \displaystyle \bigoplus_{i=1}^3 R(-d_i) &
\rightarrow & R & \rightarrow & R/J & \rightarrow & 0\\
&& \downarrow && \downarrow \\
0 & \rightarrow & F & \rightarrow & \displaystyle \bigoplus_{i=1}^4 R(-d_i) &
\rightarrow & R & \rightarrow & R/I & \rightarrow & 0\\
&& \downarrow && \downarrow \\
&& R(-d_4) & = & R(-d_4) \\ 
&&  && \downarrow \\
&&  && 0
\end{array}
\end{equation}
Note that
\[
H^1_*({\mathcal E}) \cong R/J \ \ \hbox{and} \ \ H^1_*({\mathcal F}) \cong R/I,
\]
where ${\mathcal E}$ and ${\mathcal F}$ are the sheafifications of $E$ and $F$,
respectively. Now consider the commutative diagram of locally free sheaves
\begin{equation} \label{bigdiagram2}
\begin{array}{cccccccccccccccccccc}
&& 0 && 0 && 0 \\
&& \downarrow && \downarrow && \downarrow \\
0 & \rightarrow & {\mathcal E}(-1) & \rightarrow & {\mathcal E} & \rightarrow &
{\mathcal E}_{|L} & \rightarrow & 0 \\
&& \downarrow && \downarrow && \downarrow \\
0 & \rightarrow & {\mathcal F}(-1) & \rightarrow & {\mathcal F} & \rightarrow &
{\mathcal F}_{|L} & \rightarrow & 0 \\
&& \downarrow && \downarrow && \downarrow \\
0 & \rightarrow & {\mathcal O}_{\proj{2}}(-d_4-1) & \rightarrow & {\mathcal
O}_{\proj{2}} (-d_4) & \rightarrow & {\mathcal O}_L(-d_4) & \rightarrow & 0 \\
&& \downarrow && \downarrow && \downarrow \\
&& 0 && 0 && 0 
\end{array}
\end{equation}
Twisting by $d_4 = \lambda-2$ and taking cohomology, we know by Theorem
\ref{HMNW result} that 
\[
H^1({\mathcal E}(\lambda-3)) \hookrightarrow H^1({\mathcal E}(\lambda-2))
\]
so (\ref{bigdiagram2}) becomes
\begin{equation} \label{bigdiagram3}
\begin{array}{cccccccccccccccccccc}
&& 0 && 0 && 0 \\
&& \downarrow && \downarrow && \downarrow \\
0 & \rightarrow & H^0({\mathcal E}(\lambda-3)) & \rightarrow & H^0({\mathcal
E}(\lambda-2)) & \rightarrow & H^0({\mathcal E}_{|L}(\lambda-2)) & \rightarrow &
0 \\ 
&& \downarrow && \downarrow && \phantom{\gamma} \downarrow \gamma \\
0 & \rightarrow & H^0({\mathcal F}(\lambda-3)) & \rightarrow & H^0({\mathcal
F}(\lambda-2)) & \stackrel{\beta}{\rightarrow} & H^0({\mathcal
F}_{|L}(\lambda-2)) &
\rightarrow &
\hbox{?} \\ 
&& \downarrow && \downarrow && \downarrow \\
0 & \rightarrow & 0 & \rightarrow & H^0({\mathcal O}_{\proj{2}})  &
\rightarrow & H^0({\mathcal O}_L) & \rightarrow & 0 \\ && 
 && \downarrow && \phantom{\alpha} \downarrow \alpha \\ &&  && H^1({\mathcal
E}(\lambda-2)) && H^1({\mathcal E}_{|L}(\lambda-2))
\end{array}
\end{equation}

\medskip

\noindent {\bf Claim:} The vertical map $\alpha$ in (\ref{bigdiagram3}) is an
injection.

\medskip

We will prove this claim shortly, but first we note that this completes the
proof of our desired injection, since it means that the vertical map $\gamma$
is an isomorphism, and so $\beta$ must be surjective, proving the injectivity
of $H^1({\mathcal F}(\lambda-3)) \rightarrow H^1({\mathcal F}(\lambda-2))$,
which is the desired one.

Because of Cases 1 and 2 above, we may safely assume that $d_3 < d_1 +d_2+1$. 
Then by \cite{HMNW}, Corollary 2.2, when $d_1+d_2+d_3$ is even, the splitting
type of $\mathcal E$ is $a_{\mathcal E}(\ell) = (-\lambda, -\lambda)$.  When
$d_1+d_2+d_3$ is odd, the splitting type is $(-\lambda, -\lambda+1)$.  We treat
the case when $d_1+d_2+d_3$ is even, leaving the similar odd case to the
reader.  

Let $L$ be a general line in $\proj{2}$ (and we use the same notation for the
corresponding general linear form).  We know that ${\mathcal E}_{|L} \cong
{\mathcal O}_L(-\lambda)^2$.  We have to find ${\mathcal F}_{|L}$.  Consider
the exact sequence
\[
0 \rightarrow {\mathcal O}_L(-\lambda)^2 \rightarrow {\mathcal F}_{|L}
\rightarrow {\mathcal O}_L(-\lambda+2) \rightarrow 0.
\]
Twisting and taking cohomology we get $h^0({\mathcal F}_{|L}(\lambda-1)) = 2$
and $h^1({\mathcal F}_{|L}(\lambda-1)) = 0$.  Then $h^0({\mathcal
F}_{|L}(\lambda-2))$ can only be 0 or 1.  If we show that it is 0 then this
proves the injectivity of $\alpha$ as desired.  

By considering Chern classes, we see that the only possibilities for ${\mathcal
F}_{|L}$ are
\begin{equation}\label{possibilities}
{\mathcal O}_L (-\lambda)^2 \oplus {\mathcal O}_L(-\lambda+2) \ \ \hbox{ or } \
\ {\mathcal O}_L(-\lambda+1)^2 \oplus {\mathcal O}_L(-\lambda).
\end{equation}
We claim that the first of these is impossible.  Let $\bar I = (\bar
F_1,\bar F_2, \bar F_3, \bar F_4)$ be the restriction of $I$ to $\bar R :=
R/(L)$ and consider the exact sequence
\[
0 \rightarrow H^0_*({\mathcal F}_{|L}) \rightarrow \bigoplus_{i=1}^3 \bar R(-d_i)
\oplus \bar R(-\lambda+2) \rightarrow \bar I \rightarrow 0.
\]
Suppose that ${\mathcal F}_{|L}$ is the first of the sheaves given in
(\ref{possibilities}).  Because ${\mathcal E}_{|L} = {\mathcal
O}_L(-\lambda)^2$, the summand $\bar R(-\lambda+2)$ in $H^0_*({\mathcal F}_{|L})$
cannot represent a syzygy for only $\bar F_1,\bar F_2,\bar F_3$.  But then this
means that $\bar F_4$ is not a minimal generator of $\bar I$, since its degree
is precisely $\lambda-2$.  What does this say about $F_4$ itself?  Consider the
exact sequence
\[
0 \rightarrow (R/J)_{\lambda-3} \stackrel{\times L}{\longrightarrow}
(R/J)_{\lambda-2} \stackrel{\rho}{\longrightarrow} (\bar R/\bar J)_{\lambda-2}
\rightarrow 0,
\]
where $\rho$ is the restriction map, $\rho(F) = \bar F$.  The assertion that
$\bar F_4$ is not a minimal generator of $\bar I$ means that $\bar F_4 \in \bar
J$, so  $F_4 \in \ker \rho$ (viewing $F_4$ as a non-zero element of $R/J$).  But
we assumed that $F_4 \notin \ker \rho$ for some $L$, hence this is true for the
general $L$.  This contradiction completes the proof.
\end{proof}


\section{Other appearances of the Weak Lefschetz Property}

One theme of this paper is that the Weak Lefschetz property always seems to
appear in ``general'' situations.  This section gives some other instances of
this phenomenon.

\begin{proposition} \label{hypersurf sect}
Let $X \subset \proj{2}$ be any zeroscheme, with saturated ideal $I_X$.  Let $F
\in k[x_1,x_2,x_3]_d = R_d$ be a generally chosen polynomial.  Then the Artinian
ideal $I_X + (F) =: I \subset R$ has the Weak Lefschetz property.
\end{proposition}

\begin{proof}
Let $L$ be a general linear form and let $\bar R = R/(L)$.  Let $\ell$ be the
image of $L$ in $R/I$.  We have a commutative diagram
\[
\begin{array}{ccccccccccccccc}
&&&&&& 0 \\
&&&&&&\downarrow \\
&&0 &&0 &&\ker (\times \ell) \\
&& \downarrow && \downarrow && \downarrow \\
0 & \rightarrow & (R/I_X)_{t-d} & \stackrel{\times F}{\longrightarrow} &
(R/I_X)_t & \rightarrow & (R/I)_t & \rightarrow & 0 \\ 
&& \phantom{\times L} \downarrow \times L && 
\phantom{\times L} \downarrow \times L && 
\phantom{\times \ell} \downarrow \times \ell\\
 0 & \rightarrow & (R/I_X)_{t-d+1} & \stackrel{\times F}{\longrightarrow} &
(R/I_X)_{t+1} & \rightarrow & (R/I)_{t+1} & \rightarrow & 0 \\
&& \downarrow && \downarrow && \downarrow \\
  & & (\bar R/ \bar I_X)_{t-d+1} & \stackrel{\times F}{\longrightarrow} &
(\bar R/ \bar I_X)_{t+1} & \rightarrow & \coker (\times \ell) \\
&& \downarrow && \downarrow && \downarrow \\
&&0 &&0 &&0
\end{array}
\]
Note  that $\bar R / \bar I_X$ is the Artinian reduction of $R/I_X$, and hence
has the Strong Lefschetz property by \cite{HMNW} Proposition 4.4.  Also, by the
Snake Lemma we have the exact sequence
\[
0 \rightarrow \ker (\times \ell) \rightarrow (\bar R/ \bar I_X)_{t-d+1}
\stackrel{\times \bar F}{\longrightarrow} (\bar R/ \bar I_X)_{t+1} \rightarrow
\coker (\times \ell) \rightarrow 0.
\]
Thus since $F$ is general, $\times \bar F$ has maximal rank, and hence  the same
is true of the vertical map $\times \ell$.
\end{proof}

\begin{corollary}
Let $C \subset \proj{3}$ be an arithmetically Cohen-Macaulay curve and let
$\tilde F \in S_d$ be a {\em general } homogeneous polynomial of degree $d$,
where $S = k[x_0,x_1,x_2,x_3]$.  Let $Z \subset C$ be the zeroscheme cut out by
$\tilde F$, so $I_Z = I_C + (\tilde F)$ is its saturated homogeneous ideal. 
Then {\em any} Artinian reduction of $S/I_Z$ has the Weak Lefschetz property.
\end{corollary}

\begin{proof}
If $A$ is the Artinian reduction of $S/I_Z$, we have $A \cong S/(I_Z +(L)) =
S/(I_C + (\tilde F,L))$, where $L$ is a linear form not vanishing at any point
in the support of $Z$.  Let $X$ be the hyperplane section of $C$ cut out by
$L$.  Since $L$ avoids the points of $Z$, we have that $X$ is also a
zeroscheme.  So $X \subset \proj{2} = H_L$.  Let $R = S/(L)$ and let $F$ be the
restriction of $\tilde F$ to $R$.  Note that $I_X = I_C + (L)$ and $A \cong
R/(I_X+(F))$.  Hence the result follows from Proposition \ref{hypersurf sect}.
\end{proof}

If the degree is large enough, we can improve on Proposition \ref{hypersurf
sect} by removing the assumption that $F$ be general.

\begin{proposition}
Let $X \subset \proj{2}$ be a zeroscheme with saturated ideal $I_X$ and minimal
free resolution
\[
0 \rightarrow \bigoplus_{i=1}^{r-1} R(-a_i) \rightarrow \bigoplus_{i=1}^r
R(-d_i) \rightarrow I_X \rightarrow 0.
\]
Let $F \in R_d$ be {\em any} homogeneous polynomial which does not vanish at
any point in the support of $X$.  Let $a = \max \{ a_i \}$.  If $d \geq a-1$
then $R/(I_X + (F))$ has the Weak Lefschetz property.
\end{proposition}

\begin{proof}
Suppose that $\deg X = e$.  The Hilbert function of $R/I_X$ satisfies
\[
h_{R/I_X } (t) = 
\left \{
\begin{array}{l}
\hbox{strictly increasing until $t = a-2$} \\ \\
e \hbox{ for all $t \geq a-2$}
\end{array}
\right.
\]
The Hilbert function of $R/(I_X + (F))$ is 
\[
h_{R/(I_X +(F))}(t) = h_{R/I_X}(t) - h_{R/I_X}(t-d).
\]
In particular, since we have chosen $d \geq a-1$, we have for $t \leq a-3$ that
\[
\begin{array}{ccccccc}
[R/(I_X +(F))]_t & \stackrel{\times L}{\hookrightarrow} & [R/(I_X +(F))]_{t+1} \\
|| && || \\
(R/I_X)_t && (R/I_X)_{t+1}
\end{array}
\]
For $t \geq a-2$ we have 
\[
\begin{array}{ccccccccccc}
(R/I_X)_t & \stackrel{\stackrel{\times L}{\sim}}{\longrightarrow} &
(R/I_X)_{t+1} \\
\downarrow && \downarrow \\

[R/ (I_X +(F))]_t & \longrightarrow & [ R/(I_X +(F)) ]_{t+1} \\
\downarrow && \downarrow \\
0 && 0
\end{array}
\]
which implies the desired surjectivity.
\end{proof}

\end{document}